\documentclass[a4paper,10pt]{article}
\pdfoutput=1
\makeatletter

\@addtoreset{equation}{section}
\makeatother

\usepackage[latin1]{inputenc}
\usepackage{amsfonts}
\usepackage{amsmath}
\usepackage{amssymb}
\usepackage{eufrak}
\usepackage{graphicx}
\usepackage{mathrsfs}
\usepackage{tikz}
\usepackage{calc}

\newcommand{\ds}{\displaystyle}

\newcommand{\R}{\mathbb{R}}

\newcommand{\Sc}{\mathbb{S}}

\newcommand{\be} {\begin{equation}}
\newcommand{\ee} {\end{equation}}
\newcommand{\bea} {\begin{eqnarray}}
\newcommand{\eea} {\end{eqnarray}}
\newcommand{\Bea} {\begin{eqnarray*}}
\newcommand{\Eea} {\end{eqnarray*}}

\newcommand{\al} {\alpha}

\newcommand{\de} {\delta}
\newcommand{\ga} {\gamma}

\newcommand{\om} {\omega}
\newcommand{\De} {\Delta}
\newcommand{\la} {\lambda}
\newcommand{\si} {\sigma}

\newcommand{\La} {\Lambda}

\newcommand{\va} {\varphi}
\newcommand{\eps} {\varepsilon}

\DeclareMathOperator{\ud} {d \!}
\renewcommand{\le}{\leqslant}
\renewcommand{\ge}{\geqslant}

\newcommand{\GL}{{\bf GL}}

\newcommand{\Or}{{\bf O}}

\newcommand{\Mnm} {\mathcal{M}_{n,m}}
\DeclareMathOperator{\tr} {Tr}

\newcommand{\LL}{\mathscr{L}^2}
\newcommand{\reg}{\mathrm{Reg}}
\newcommand{\m}{\underline{m}}
\newcommand{\vdra}{\mathscr{D}_d}
\newcommand{\sD}{s_{\la}}
\newcommand{\sL}{\mathscr{S}_\la(L)}

\newcommand{\pd}{{\begin{tikzpicture} 
\draw (0,0)--+(3pt,0)--+(3pt,3pt)--+(0pt,3pt)--cycle ;
\draw (3pt,0)--+(3pt,0)--+(3pt,3pt)--+(0pt,3pt)--cycle;
\end{tikzpicture}}}

\newcommand{\pq}{{\begin{tikzpicture} 
\draw (0,0)--+(3pt,0)--+(3pt,3pt)--+(0pt,3pt)--cycle ;
\draw (3pt,0)--+(3pt,0)--+(3pt,3pt)--+(0pt,3pt)--cycle;
\draw (6pt,0)--+(3pt,0)--+(3pt,3pt)--+(0pt,3pt)--cycle;
\draw (9pt,0)--+(3pt,0)--+(3pt,3pt)--+(0pt,3pt)--cycle;
\end{tikzpicture}}}

\newcommand{\pdd}{{\begin{tikzpicture} 
\draw (0,0)--+(3pt,0)--+(3pt,3pt)--+(0pt,3pt)--cycle ;
\draw (3pt,0)--+(3pt,0)--+(3pt,3pt)--+(0pt,3pt)--cycle;
\draw (0pt,3pt)--+(3pt,0)--+(3pt,3pt)--+(0pt,3pt)--cycle;
\draw (3pt,3pt)--+(3pt,0)--+(3pt,3pt)--+(0pt,3pt)--cycle;
\end{tikzpicture}}}

\newcommand{\pdi}[2]{{\raisebox{-5pt}{\begin{tikzpicture} 
\draw (0,0) rectangle node [anchor = mid] {$#1$} +(12pt,12pt)  ;
\draw (12pt,0) rectangle node [anchor = mid] {$#2$} +(12pt,12pt) ;
\end{tikzpicture}}}}

\newcommand{\pqi}[4]{{\raisebox{-5pt}{\begin{tikzpicture} 
\draw (0,0) rectangle node [anchor = mid] {$#1$} +(12pt,12pt)  ;
\draw (12pt,0) rectangle node [anchor = mid] {$#2$} +(12pt,12pt) ;
\draw (24pt,0) rectangle node [anchor = mid] {$#3$} +(12pt,12pt)  ;
\draw (36pt,0) rectangle node [anchor = mid] {$#4$} +(12pt,12pt) ;
\end{tikzpicture}}}}

\newcommand{\pddi}[4]{{\raisebox{-10pt}{\begin{tikzpicture} 
\draw (0,0) rectangle node [anchor = mid] {$#1$} +(12pt,12pt)  ;
\draw (12pt,0) rectangle node [anchor = mid] {$#2$} +(12pt,12pt) ;
\draw (0pt,12pt) rectangle node [anchor = mid] {$#3$} +(12pt,12pt)  ;
\draw (12pt,12pt) rectangle node [anchor = mid] {$#4$} +(12pt,12pt) ;
\end{tikzpicture}}}}

\newcommand{\pdui}[3]{{\raisebox{-9.5pt}{\begin{tikzpicture} 
\draw (0,0) rectangle node  {$#1$} +(12pt,12pt)  ;
\draw (12pt,0) rectangle node  {$#2$} +(12pt,12pt) ;
\draw (0pt,12pt) rectangle node  {$#3$} +(12pt,12pt) ;
\end{tikzpicture}}}}

\newcommand{\bitab}[2]{\left\{ #1 ,  #2\right\}} 
\newcommand{\psp}[2]{\left[ #1 ,  #2\right]} 
\newcommand{\psll}[2]{\left\langle #1 ,  #2\right\rangle} 

\usepackage{theorem}
\newtheorem{thm}{Theorem}[section]
\newtheorem{prop}[thm]{Proposition}
\newtheorem{lm}[thm]{Lemma}
\newtheorem{cor}[thm]{Corolary}
\theorembodyfont{\upshape}
\newtheorem{defi}[thm]{Definition}
\newtheorem{ex}[thm]{Example}

\theorembodyfont{\sffamily}

\newenvironment{dem}{\noindent  {\sl Demonstration.{\space \it }}}{\hfill $\square$}
\newenvironment{proof}{\noindent {\sl Proof.{\space \it }}}{\hfill $\square$}

\title{Extreme lattices and vexillar designs\footnote{Version of \today. MSC2000 : 05B30, 14M15, 20C15.}}
\author{Bertrand Meyer\thanks{\texttt{bertrand.meyer@math.u-bordeaux1.fr} -- Institut de mathématiques de Bordeaux -- UMR~5251 -- Université Bordeaux 1 -- 351, cours de la libération -- 33 405 Talence Cedex -- France}}
\date{}

\begin{document}
\maketitle

\begin{abstract} We define a notion of vexillar design for the flag variety in the spirit of the spherical designs introduced by Delsarte, Goethals and Seidel. For a finite subgroup of the orthogonal group, we explain how conditions on the group have the orbits of any flag under the group action be a design and point out why the minima of a lattice in the sense of the general Hermite constant forming a 4-design implies being extreme. The reasoning proves useful to show the extremality of many new expected examples ($E_8$, $\La_{24}$,  Barnes-Wall lattices, Thompson-Smith lattice for instance) that were out of reach until now.
\end{abstract}

A very general form of Hermite constant associated with an algebraic group over some number field and a strongly rational irreducible representation was introduced in \cite{W}. The framework of this constant is large enough to encompass all the previously studied generalisations of the Hermite constant, for instance the Rankin constant or the Humbert constant. Futhermore, this new definition provides a good point of view to tackle number theoretic issues such as looking for rational points of minimal height on the flag variety or writing results in the spirit of Siegel lemmas. Yet to that day, very few is known about explicit values of these constants.

In this article, we shall stay on the field of rational numbers and escape any number theoretical background. We choose once for all an integer $n$ and some non-increasing integers $(\breve{\la}_i)_{1 \le i \le \breve{s}}$ less than~$n$. The Hermite constant associated with the general linear group and the polynomial representation of weight $\la$ can be expressed after reformulation as follows. For any full-rank lattice~$L$ contained in~$\R^n$, we define
\be \ga(L) = \inf_{\La \subset L} \frac{\det(\La_1) \cdots \det(\La_{\breve{s}}) }{\left( \det (L) \right) ^{\frac{|\la|}{n}}}\ee
where $\La = (\La_1, \dots, \La_{\breve{s}})$ is a chain of nested sublattices of $L$ satisfying the condition $\text{rk} (\La_i) = \breve{\la}_i$ for all $1 \le i \le \breve{s}$. The Hermite constant we are interested in is simply 
\be \ga_{\la,n} = \sup_{L} \ga(L) \ee
where the supremum is  taken over all the full-rank lattices $L$ of $\R^n$.

A lattice $L$ is called extreme when it achieves a local maximum of $\ga(L)$. A complete theory to characterise the extreme forms with respect to $\la$ has been issued in \cite{Mey}. Yet, to that day, very few examples of extreme lattices have been displayed in the general case : one of the goal of this paper is to prove that lattices such as for instance the Leech lattice, some of the root lattices, the Barnes-Wall lattices or the the Thompson--Smith lattice are indeed extreme for any weight $\la$. Considering the ubiquity in the literature of these examples, they provide also good lower bounds for the general Hermite constant and serious candidates to be actually the global maxima of~$\ga(L)$.

 A notion of spherical design was introduced by Delsarte, Goethals and Seidel in the paper~\cite{DGS}. Later, Boris Venkov set forth in a milestone article~\cite{V} the link between spherical designs and some special lattices he calls strongly perfect. Lempken, Schröder and Tiep~\cite{LST} explored designs arising as orbits of finite groups and used the automorphism group of lattices to conclude via designs that they are extreme in the usual sense. This procedure was further used in~\cite{BCN} for Grassmanians in relation with the Rankin invariant and the method of Delsarte, Goethals and Seidel was further expanded in~\cite{BBC}. Besides the groups that provide Grassmannian $4$-designs by orbits and thus extreme lattices whenever they occur as the automorphism group of the lattice were classified in \cite{TIE}.

In section \ref{s:prlg}, we recall some notions on the flag variety and how the space of regular functions on it decomposes into orthogonal irreducible components under the action of the orthogonal group. This can be explicitly performed by using determinantal monomials and Young tableaux. Then in section \ref{s:zonfun}, we construct and explicit some zonal functions of low degree. This enables us to define in section \ref{s:desflag} a notion of vexillar design, which suits our needs for the general Hermite constant. In particular, we show some equivalent conditions for a finite set to be a design. It makes it possible then to tie the absence of invariants of a finite group to an orbit under the group being a design. In section \ref{s:sperf}, we exhibit the link between $4$-designs and extreme lattices. In the last section, we describe some examples which benefits from this theory and show for instance that the root lattice $E_8$, the Leech lattice $\La_{24}$ or the Barnes-Wall lattices or the Thompson--Smith lattice are new extreme lattices.

\section{Prolegomena} 
\label{s:prlg}

Let us call \emph{flag} of the vector space~$\R^n$ of shape~$d=(d_1, \dots, d_\ell)$ and let us denote~$\De$ any sequence of embeded subspaces
\be \De :  \quad \{0\} \subsetneq V_\ell \subsetneq \cdots \subsetneq V_i \subsetneq \cdots \subsetneq V_1 \subsetneq \R^n, \ee
the dimension $d_i = \dim V_i$ of which has been fixed once for all.  Let us take  $m=d_1$, $m_1 = d_\ell$, $m_2 = d_{\ell-1}-d_\ell$,\dots  $m_i = d_{\ell+1-i}-d_{\ell+2-i}$ when $i$ lies between $2$ and $\ell$ and $m_{\ell+1} = n-d_1$.  The set of all flags of shape~$d$ constitutes a variety, denoted~$\vdra$ ; it identifies with the quotient~$\Or(n)/ \Or(m_1)\times \cdots \times \Or(m_{\ell+1})$.  We shall mark in the sequel a flag~$\De$ by a matrix with orthonormal columns~$X_\De \in \Mnm(\R)$ the $d_i$ first column vectors of which span the space $V_i$ ; the matrix $X_\De$ is not unique but these matrices deduce one from the other by a right multiplication of $\Or(\m) = \Or(m_1) \times \cdots \Or(m_\ell)$.

Let us call \emph{partition} any non-increasing finite sequence of natural integer $\mu = (\mu_1, \dots, \mu_u)$. The \emph{degree}, denoted~$|\mu|$, of a partition is the sum of its parts~; the \emph{depth}, denoted~$\lfloor \mu \rfloor$, is the number of its parts. The transpose partition of~$\mu$ is denoted $\breve{\mu} = (\breve{\mu}_1, \dots , \breve{\mu}_{\breve{u}})$. A partition can be represented by its Ferrer diagram, which is a drawing in the first quadrant of $u$ rows of boxes of lengths $(\mu_i)_{1\le i \le u}$ or else of $\breve{u}$ columns of lenghts $(\breve{\mu}_i)_{1\le i \le \breve{u}}$. A \emph{tableau}~$T$ is the data of a Ferrer diagram inscribed with natural integers, $T(i,j)$ referring to the entry located at the abscissa~$i$ and the ordinate~$j$. A tableau is said \emph{standard} when the integers written in the diagram are increasing along the columns and non-decreasing along the rows. The content~$C_T$ of a tableau~$T$ is the count for any integer of its number of occurence in~$T$. Eventually, we call \emph{bitableau}~$B$ any pair of tableaux~$\bitab{T}{\Theta}$. In this article, we shall only consider bitableaux with entries between  $1$ and $n$ on the left and between $1$ and $m$ on the right. A bitableau is said standard if the two tableaux that compound it are standard. The content of a bitableau~$B$ is the pair of its contents~$(C_T,C_\Theta)$.

The letter $X$ will always indicate the matrix
$$ X = \begin{pmatrix} x_{1,1} & \cdots & x_{1,m} \\ \vdots & & \vdots \\ x_{n,1} & \cdots & x_{n,m} \end{pmatrix} $$
the entries $(x_{i,j})_{1 \le i \le n, 1 \le j \le m}$ of which are unknowns and $\reg(\Mnm(\R))$ the space of \emph{regular functions} in the $n\cdot m$ variables. This space is equiped with the \emph{Euclidean scalar product}~$\psp{\cdot}{\cdot} $ that makes orthogonal any two distincts monomials and for which the scalar square of a monomial~$x_{1,1}^{\al_{1,1}} \cdots x_{n,m}^{\al_{n,m}}$ is equal to $\al ! / |\al|!$, where $\al ! = \al_{1,1} ! \cdots \al_{n,m} !$. Of a bitableau~$\bitab{T}{\Theta}$ we can form the following \emph{determinantal monomial}, said of shape~$\mu$,
\be \mathscr{M}_{\bitab{T}{\Theta}} = \prod_{i=1}^{\breve{u}} \det \begin{pmatrix} x_{T(i,j), \, \Theta(i,j')} \end{pmatrix}_{1 \le j, \, j' \le \breve{\mu}_i} \ee
which is a product of minors of $X$, where the choice of the excerpted rows is commanded by the a column from the left tableau $T$ and the choice of the columns by the concomitant column from the right tableau $\Theta$. As a consequence of the straightening law by Doubilet {\it et al.} \cite{DRS}, the set of all the standard monomials, {\it i. e.} the monomials $\mathscr{M}_B$ where $B$ is a standard bitableau, sets up a basis of the Hilbert space of regular functions $\reg(\Mnm(\R))$.

We denote by~$I_t(X)$ the ideal of $\reg(\Mnm(\R))$ spanned by the minors of size~$t$ of $X$. For a partition $\mu$, $I^{(\mu)}(X)$ denotes the product of ideals $I^{(\mu)} (X) = I_{\breve{\mu}_1}(X) \cdots I_{\breve{\mu}_s}(X)$. A criterion makes it easy to check the membership to this ideal. Let us call, for any partition~$\mu$, $\varkappa_t(\mu) = \sum_{k=t}^{+\infty} {\mu}_k$ the number of boxes of~$\mu$ above the $t$\textsuperscript{th} column. A monomial~$\mathscr{M}$ of shape~$\si$ belongs to the ideal~$I^{(\mu)}(X)$ if and only if, for any integer~$t$, the inequality $\varkappa_t(\si) \ge \varkappa_t(\mu)$ holds (see \cite{DCEP}, read also chap. 11 of \cite{BV}). We shall also denote~ $I^{(\mu)}_>(X)$ the ideal spanned by the ideals~$I^{(\si)}(X)$ indexed by the partitions $\si$ which strictly contain~$\mu$ in the sense of the inclusion of Ferrer diagram.

We recall that the \emph{polynomial representations} of $\GL_n(\R)$ are parametrised by the set of all \emph{partitions}~$\mu$ with less than $n$ parts (see \cite{Fu} for more details). We denote by~$\Sc^\mu$ the \emph{Schur functor} : the space~$\Sc^\mu(\R^n)$  is one of the irreducible polynomial representations of the group~$\GL_n(\R)$ ; all of them appear like that ; $\mu$~is called the weight of the representation. The space~$\Sc^\mu(\R^n)$ can be represented by formal linear combinations of diagrams of shape~$\mu$ inscribed with vectors, on which directly apply the elements of~$\GL_n(\R)$ when the group acts on the space. The representatives are not unique in general. We may refer to the elements as flag vectors if there is only one term in the combination. If a basis~$(e_i)_{1 \le i \le n}$ of~$\R^n$ is decided, the flag vectors~$e_T$ of $\Sc^\mu(\R^n)$, inscribed in position~$(i,j)$ of the vector~$e_{T(i,j)}$, provide in particular a basis of~$\Sc^\mu(\R^n)$ when $T$ runs through the set of standard Young tableaux with coefficients between $1$ and $n$. In the same way, if a basis~$(\eps_i)_{1 \le i \le m}$ of $\R^m$ is fixed, we are afforded a basis by flag vectors~$e_\Theta$ of the space~$\Sc^\mu(\R^m)$ when $\Theta$ runs among the standard Young tableaux with coefficients between $1$ and $m$.  We recall also that the Schur module $\Sc^\mu(\R^n)$ is simply the symmetric power $\mathsf{Sym}^{|\mu|}(\R^n)$ when $\mu$ consists of only one part and the alternating power $\mathsf{Alt}^{|\mu|}(\R^n)$ when $\breve{\mu}$ consists of a unique part.

The space of regular functions on $\Mnm(\R) \cong \R^n \otimes \R^m$ naturally possesses a structure of $\GL_n(\R) \times \GL_m(\R)$-module under the action of the left multiplication of the variable by an element of~$\GL_n(\R)$ and right multiplication by an element of $\GL_m(\R)$. The ideals $I^{(\mu)}(X)$ are invariant under this action ; the scalar product~$\psp{ \cdot}{ \cdot} $ is invariant under the action of the subgroup~$\Or(n)\times \Or(m)$.  The decomposition of this space of functions into irreducible $\GL_n(\R) \times \GL_m(\R)$-modules is well known \cite{BV} and can be state as follows : 
\be \reg(\Mnm(\R) ) = \bigoplus_{\mu}^{\perp} M^\mu \cong \bigoplus_{\substack{\mu}}^\perp \Sc^{\mu} (\R^n) \otimes \Sc^\mu (\R^m), \ee where $\mu$ describes the partitions with less than $m$ and $n$ parts and the isomorphism has to be understood componentwise. It can be further demonstrated that~$M^\mu$ is the unique complement invariant under $\GL_n(\R) \times \GL_m(\R)$ of the ideal~$I^{(\mu)}_>(X)$ in~$I^{(\mu)}(X)$ (see sect. 3 of \cite{DCEP}). The vector~$e_{T_0} \otimes e_{\Theta_0}$ of~$\Sc^\la(\R^n) \otimes \Sc^\la(\R^m)$ identifies to the polynomial $\mathscr{M}_{\bitab{T_0}{\Theta_0}}$ of $\reg(\Mnm(\R))$. Thus, the isomorphism between $\Sc^{\mu} (\R^n) \otimes \Sc^\mu (\R^m) $ and $M^\mu$ can by totally explicitly written by letting act on both sides the groups~$\GL_n(\R)$ and~$\GL_m(\R)$. As a consequence of the straightening law by Doubilet {\it et al.} \cite{DRS}, we have even 
$$ e_T \otimes e_\Theta  \mapsto  \mathscr{M}_{\bitab{T}{\Theta}}   + I^\mu_>(X) $$ 
for any pair of flag vectors $e_T$ and $e_\Theta$. Nevertheless in general, the equality holds only modulo $ I^\mu_>(X)$ and $\mathscr{M}_{\bitab{T}{\Theta}}$ does not necessarily belong to $M^\mu$, as shows the following example. Let us take 
$$T_1 = \pdui{1}{3}{2}, \quad 
T_2 = \pdui{1}{2}{3}, \quad
T_3 = \pdui{3}{1}{2}, \quad 
\Theta_1 = \pdui{1}{2}{3}, \quad
 B^\sharp = {\bitab{ 
\raisebox{-15pt}{ {\begin{tikzpicture} 
\draw (0,0) rectangle node  {1} +(12pt,12pt)  ;
\draw (0pt,12pt) rectangle node  {2} +(12pt,12pt) ;
\draw (0pt,24pt) rectangle node  {3} +(12pt,12pt) ;
\end{tikzpicture}}} }{ 
\raisebox{-15pt}{\begin{tikzpicture} 
\draw (0,0) rectangle node  {1} +(12pt,12pt)  ;
\draw (0pt,12pt) rectangle node  {2} +(12pt,12pt) ;
\draw (0pt,24pt) rectangle node {3} +(12pt,12pt) ;
\end{tikzpicture}} }} .$$
Note that $T_1$, $T_2$ and $\Theta_1$ are standard but not $T_3$. Then
$$\mathscr{M}_{\bitab{T_3}{\Theta_1}} =
\mathscr{M}_{\bitab{T_1}{\Theta_1}} 
-
\mathscr{M}_{\bitab{T_2}{\Theta_1}}  
+ \mathscr{M}_{B^\sharp}
$$
and $\mathscr{M}_{B^\sharp} \in I^{{\begin{tikzpicture} 
\draw (0,0)--+(3pt,0)--+(3pt,3pt)--+(0pt,3pt)--cycle ;
\draw (0pt,3pt)--+(3pt,0)--+(3pt,3pt)--+(0pt,3pt)--cycle;
\draw (0pt,6pt)--+(3pt,0)--+(3pt,3pt)--+(0pt,3pt)--cycle;
\end{tikzpicture}}} (X) \subseteq I_{>}^{{\begin{tikzpicture} 
\draw (0,0)--+(3pt,0)--+(3pt,3pt)--+(0pt,3pt)--cycle ;
\draw (0pt,3pt)--+(3pt,0)--+(3pt,3pt)--+(0pt,3pt)--cycle;
\draw (3pt,0)--+(3pt,0)--+(3pt,3pt)--+(0pt,3pt)--cycle;
\end{tikzpicture}}} (X)$.

We shall denote $\phi(e_T \otimes e_\Theta)$ the unique representative of the class $\mathscr{M}_{\bitab{T}{\Theta}} + I^{(\mu)}_>(X)$ in the irreducible $\GL_n (\R) \times \GL_m (\R) $-module $M^\mu$.

Let us mention that when in one of the two sides of the bitableau all the indices of any column collectively spread among all the longer columns, then the monomial~$\mathscr{M}_{\bitab{T}{\Theta}}$ really belongs to $M^\mu$. This remark makes it possible to recover one of the classical construction of~$\Sc^{\mu}(\R^n)$ : by systematically choosing for $\Theta$ the so-called \emph{initial} tableau
$$\Theta_0 = \raisebox{-0.5\height + 3pt}{\begin{tikzpicture} 
\draw (0,0) rectangle node {$1$} +(12pt,12pt)  ;
\draw (12pt,0) rectangle node  {1} +(12pt,12pt) ;
\draw (24pt,0) rectangle node  {1} +(12pt,12pt) ;
\draw (36pt,0) rectangle node  {$\cdots$} +(12pt,12pt) ;
\draw (48pt,0) rectangle node  {1} +(12pt,12pt) ;
\draw (0pt,12pt) rectangle node  {$2$} +(12pt,12pt)  ;
\draw (12pt,12pt) rectangle node  {2} +(12pt,12pt) ;
\draw (24pt,12pt) rectangle node  {$\cdots$} +(12pt,12pt) ;
\draw (36pt,12pt) rectangle node  {2} +(12pt,12pt) ;
\draw (0pt,24pt) rectangle  node {\! :} +(12pt,12pt)  ;
\draw (12pt,24pt) rectangle  +(12pt,12pt) ;
\draw (0pt,36pt) rectangle node  {$t$} +(12pt,12pt)  ;
\end{tikzpicture}},
$$
the monomials $\mathscr{M}_{\bitab{T}{\Theta_0}}$ span a space isomorphic to $\Sc^{\mu}(\R^n)$.

Let us observe in particular that the component of homogeneous polynomials of degree~$k$ of~$\reg(\Mnm(\R))$ exactly corresponds to the sum of the spaces indexed by a partition~$\mu$ of degree~$|\mu|= k$. 

Let $\LL(\vdra)$ be the space of \emph{square integrable functions} on the flag variety~$\vdra$ equiped with the Haar measure of total measure 1. Any regular function of~$\reg(\Mnm(\R))$ stable under the action on the right by~$\Or(\m)$ affords by restriction a function of~$\LL(\vdra)$. The space $\LL(\vdra)$ is equiped with the scalar product $\psll{f}{g} = \int_{\vdra} f\cdot g$, which is proportional to the scalar product $\psp{ \cdot}{\cdot }$ on any irreducible sub-$\Or(n) \times \Or(\m)$-module of~$\LL(\vdra)$ .
 
As an $\Or(n)$-module, the space $\Sc^{\mu}(\R^n)$ decomposes into a sum
$$ \Sc^{\mu}(\R^n) = \Sc^{[\mu]}(\R) \oplus J^\la, $$
where $\Sc^{[\la]}$ is the irreducible representation of weight $\mu$ of $\Or(n)$ if $\mu$ has less than $n/2$ parts and vanishes otherwise and $J^{\la}$ is a sum of irreducible representations of $\Or(n)$ the weight of which are strictly less than  $\mu$ (using $\varkappa_t$ in the same sense as above). Moreover, $J^{\la}$ is spanned (see section 8 in \cite{LT}) by the set of the sums~$V$ in~$\Sc^\mu(\R^n)$ of tableaux inscribed with vectors constructed as follows :  from a tableau vector~$V_0$ of $\Sc^\mu(\R^n)$,
we select $r$ boxes in each of two distinct columns, then $V$ is the sum of the tableaux $V_0$ in which for any multi-index $(i_1 < \cdots < i_r)$ of integers taken between $1$ and $n$ has been substituted by respecting the order the content of the selected boxes by the sequence of vectors~$e_{i_1}$,\,\dots \,, $e_{i_r}$ of the basis of~$\R^n$.

Let $e_\Theta$ be a basis vector of $\Sc^\mu(\R^m)$ for some standard Young tableau $\Theta$ and let $C^r_1$ and $C^r_2$ be the $r$ selected boxes in two different columns of $V_0$.  It is clear from the description of  $J^\mu$ that the polynomial $\phi( V\otimes e_\Theta)$ comprises the following factor
$$\sum_{I =  i_1 < \cdots < i_r }  \det \left( \big(x_{i,j}\big)_{\substack{i \in I \\j \in C^r_{1} }}\right) \det\left( \big(x_{i,j}\big)_{\substack{i \in I \\j \in C^r_{2} }} \right) \mod I^{(\mu)}_>$$
which is quite simply the expression of the scalar product in $\bigwedge^r \R^n$ of $\bigwedge_{j \in C^r_{1}} u_j$ and $\bigwedge_{j \in C^r_{2}} u_j$ where $u_j$ denotes the $j$\textsuperscript{th} column vector of $X$. This expression is invariant under the action of $\Or(n)$ (see section 5 of \cite{DCP}). Evaluated in the matrix~$X_\De$ which represents the flag $\De$, this expression equals to $1$ or $0$ depending wether $C^r_{1}$ and  $C^r_{2}$ are equal or not since the columns of $X_\De$ are orthonormalised. Thus the polynomials of $J^\mu \otimes \Sc^{\mu}(\R^m)^{\Or(\m)}$ identify as polynomial functions on $\vdra$ with polynomials of lower degree and we get the following decomposition of $\LL(\vdra)$ :
\be \LL(\vdra) \cong \bigoplus_{\mu} \Sc^{[\mu]} (\R^n) \otimes N^{\mu} \cong \bigoplus_{\mu} \Sc^{[\mu]} (\R^n) ^{\oplus n_\mu} \ee
where $\mu$ runs through the partitions with less than $m$ parts and $n_\mu$ is the dimension, possibly zero, of the invariant subspace $N^\mu = \Sc^{\mu}(\R^m)^{\Or(\m)}$. 

Let us notice that we also dispose of a basis of~$\Sc^{[\mu]}(\R^n)$. We say that $T$ contains a \emph{violation} if when the same integer~$v$ appears in two distinct columns, there are more indices less than~$v$ repeated in both columns than there are indices less than~$v$ absent of the two columns. Let us denote $[e_T]$ the orthogonal projection on $\Sc^{[\mu]}(\R^n)$ of $e_T$. The set of vectors $[e_T]$ without violation provides a basis of $\Sc^{[\mu]}(\R^n)$ (see \cite{LT}).

Let us eventually describe some of the stable subspaces~$N^\mu$. To that end, we can calculate for instance the image of the endomorphism of $\Sc^{\mu}(\R^m)$ defined by averaging on the full group $\va : V \mapsto \int_{\Or(\m)} \tau.V \, \ud\tau$. If $\Theta$ is a standard tableau, the content of which is not even for some integer $v$, then, swaping in $e_\Theta$ the vector $e_v$ for its opposite, we note that  $\va(e_T) = \va(-e_T)  = 0$. Thus
$n_{{\begin{tikzpicture} 
\draw (0,0)--+(3pt,0)--+(3pt,3pt)--+(0pt,3pt)--cycle ;
\end{tikzpicture}}}$, 
$n_{{\begin{tikzpicture} 
\draw (0,0)--+(3pt,0)--+(3pt,3pt)--+(0pt,3pt)--cycle ;
\draw (3pt,0)--+(3pt,0)--+(3pt,3pt)--+(0pt,3pt)--cycle ;
\draw (6pt,0)--+(3pt,0)--+(3pt,3pt)--+(0pt,3pt)--cycle ;
\end{tikzpicture}}}$, 
$n_{{\begin{tikzpicture} 
\draw (0,0)--+(3pt,0)--+(3pt,3pt)--+(0pt,3pt)--cycle ;
\draw (3pt,0)--+(3pt,0)--+(3pt,3pt)--+(0pt,3pt)--cycle;
\draw (0pt,3pt)--+(3pt,0)--+(3pt,3pt)--+(0pt,3pt)--cycle ;
\end{tikzpicture}}}$, 
$n_{{\begin{tikzpicture} 
\draw (0,0)--+(3pt,0)--+(3pt,3pt)--+(0pt,3pt)--cycle ;
\draw (0pt,3pt)--+(3pt,0)--+(3pt,3pt)--+(0pt,3pt)--cycle;
\draw (0pt,6pt)--+(3pt,0)--+(3pt,3pt)--+(0pt,3pt)--cycle ;
\end{tikzpicture}}}$, 
$n_{{\begin{tikzpicture} 
\draw (0,0)--+(3pt,0)--+(3pt,3pt)--+(0pt,3pt)--cycle ;
\draw (3pt,0)--+(3pt,0)--+(3pt,3pt)--+(0pt,3pt)--cycle ;
\draw (0pt,3pt)--+(3pt,0)--+(3pt,3pt)--+(0pt,3pt)--cycle ;
\draw (0pt,6pt)--+(3pt,0)--+(3pt,3pt)--+(0pt,3pt)--cycle ;
\end{tikzpicture}}}$, 
$n_{{\begin{tikzpicture} 
\draw (0,0)--+(3pt,0)--+(3pt,3pt)--+(0pt,3pt)--cycle ;
\draw (0pt,3pt)--+(3pt,0)--+(3pt,3pt)--+(0pt,3pt)--cycle ;
\draw (0pt,6pt)--+(3pt,0)--+(3pt,3pt)--+(0pt,3pt)--cycle ;
\draw (0pt,9pt)--+(3pt,0)--+(3pt,3pt)--+(0pt,3pt)--cycle ;
\end{tikzpicture}}}$
all vanish. 
\begin{itemize}
\item  $N^\pd$ is of dimension $n_\pd = \ell$ and is spanned by the vectors $\ds  \sum_{ d_i \le v < d_{i-1} } \pdi{\eps_v}{\eps_v} $ for $i$ between $1$ and $\ell$. Indeed, the vectors $\pdi{\eps_v}{\eps_v}$ and $\pdi{\eps_{v'}}{\eps_{v'}}$ are in the same orbit under $\Or(\m)$ if $d_i \le v,v' < d_{i-1} $, thus $\va(\pdi{\eps_v}{\eps_v}) = \va(\pdi{\eps_{v'}}{ \eps_{v'}})$. So, $\ds \va (\pdi{\eps_v}{\eps_v}) =  \frac{1}{d_{i-1}-d_i}  \sum_{d_i \le  v < d_{i-1} }
\pdi{\eps_v}{\eps_v}$. But $\ds  \sum_{d_i \le  v < d_{i-1}} \pdi{\eps_v}{\eps_v}$ is invariant under $\Or(\m)$, which completes the justification. 
\item $N^\pq$ is of dimension $n_\pq=\dfrac{\ell \, (\ell+1)}{2}$ and is spanned by the vectors $\ds \sum_{\substack{d_i \le  v < d_{i-1}  \\ d_i \le w < d_{j-1} }}
\pqi{\eps_v}{\eps_v}{\eps_w}{\eps_w} $ for $i\le j$ between $1$ and $\ell$. 
\item $N^\pdd$ is of dimension $n_\pdd = \ell$ and is spanned by the vectors $\ds  \sum_{d_i \le v < w < d_{i-1}} \pddi{\eps_v}{\eps_v}{\eps_w}{\eps_w}$ for $i\le j$ between $1$ and $\ell$.
\end{itemize}
Finally, the space of polynomial functions of degree less than 4 is reduced to
$$  \R\;  \oplus \;  \Sc^{[\pd]}(\R^n) \otimes N^\pd \; \oplus \; \Sc^{[\pdd]}(\R^n) \otimes N^\pdd \; \oplus \; \Sc^{[\pq]}(\R^n) \otimes N^\pq .$$
We denote for the rest of the text $\Upsilon_{d_{\iota}} =\ds  \sum_{v \le d_{\iota} } \pdi{\eps_v}{\eps_v} \in N^\pd $ for $1 \le \iota \le s$ ; these vectors form a basis of $N^\pd$.

In the sequel, we shall always suppose that $\mu$ possesses less than $m/2$ parts.

\section{Zonal functions}
\label{s:zonfun}

\begin{defi} For any non-zero sub-$\Or(n)$-module $\mathscr{V}$ of $\LL(\vdra)$, we call \emph{zonal function} a function $Z : \vdra \times \vdra \to \R\ $ which satisfies that
\begin{enumerate}
\item for any flag $\De_0\in \vdra$, $Z(\De_0, \cdot)$ and $Z(\cdot, \De_0)$ belong to $\mathscr{V}$,
\item for any orthogonal transformation $\tau \in \Or(n)$, we have $\; Z(\tau \De, \tau \De') = Z(\De, \De') $.
\end{enumerate}
We denote $\mathcal{Z}(\mathscr{V})$ their set.
\end{defi}

The construction below give proof of the non-emptyness of this space.

Let $(h_i)_{1 \le i \le N}$ be an orthonormal basis of $\Sc^{[\mu]}(\R^n)$, in the sense of the pull-back scalar product of $\phi(\Sc^\mu(\R^n) \otimes e_{\Theta_0})$. 
Given $\Xi$ and $\Xi'$ in $N^\mu$, we define the map~$Z_{\Xi, \Xi'}(\De, \De') $ by
\be  Z_{\Xi, \Xi'}(\De, \De')  =  \frac{1}{N} \sum_{i=1}^N \phi(h_i \otimes \Xi) (\De) \cdot  \phi(h_i \otimes \Xi' ) (\De')  .\ee
The expression~$Z_{\Xi, \Xi'}$ is bilinear in $\Xi$ and $\Xi'$.

\begin{prop} \label{prop:indep} Let $\Xi$ and $\Xi'$ be two elements of $N^\mu$,
\begin{enumerate}
\item $Z_{\Xi, \Xi'}(\De, \De')$ does not depend on the choice of the orthonormal basis $(h_i)_{1 \le i \le N}$.  
\item The application $Z_{\Xi, \Xi'}$ belongs to $\mathcal{Z}(\LL(\vdra)) $ and  bijectively sends $\phi(\Sc^{[\mu]} \otimes \Xi )$ onto $\phi(\Sc^{[\mu]}\otimes \Xi')$ by convolution.
\end{enumerate}
\end{prop}

\begin{dem} Let $(h_i')_{1 \le i \le N}$ be an other orthonormal basis of $\Sc^{[\mu]}(\R^n)$, which we express in terms of the first one by  $h_i' = \sum_{j = 1}^N \al_{j,i} h_i$ for all $i$. Then, by linearity of the morphism $\phi$ and by orthogonality of the coefficients~$((\al_{j,i}))_{1 \le j,i \le N}$,
\Bea && \sum_{i=1}^N \phi(h_i' \otimes \Xi) (\De) \cdot  \phi(h_i' \otimes \Xi' ) (\De') \\
 &=& \sum_{1 \le i,j,j' \le N} \al_{j,i} \al_{j',i} \phi(h_j \otimes \Xi) (\De) \cdot  \phi(h_{j'} \otimes \Xi' ) (\De')\\
&=& \sum_{1 \le j,j' \le N}\left( \sum_{i=1}^N \al_{j,i} \al_{j',i} \right) \phi(h_j \otimes \Xi) (\De) \cdot  \phi(h_{j'} \otimes \Xi' ) (\De')\\
&=& Z_{\Xi, \Xi'}(\De,\De'). \Eea

The second point is nothing more than a rephrasing of the first, taking into account that the scalar product is $\Or(n)$-invariant.
\end{dem}

\begin{lm} Let $\mathscr{V}$ be an sub-$\Or(n)$-module of $\LL(\vdra)$. If the dimension of $\mathcal{Z}(\mathscr{V})$ is less than 1, then $\mathscr{V}$ is irreducible.
\end{lm}

\begin{proof} This lemma is quite classic (see \cite{V}). As a consequence, $\mathcal{Z}(\LL(\vdra))$ is fully described by the functions $(Z_{\Xi, \Xi'})_{\Xi, \Xi' \in N^\mu}$ when $\mu$ varies. \end{proof}

\paragraph{Calculations}

In the case of the Grassmanians, {\it i. e.} when $\ell =1$, the multiplicities of the isotypic spaces never exceed 1 ; the zonal functions have been calculated in \cite{JC}. We compute some zonal functions for $\mu = {\begin{tikzpicture} \draw (0,0)--+(5pt,0)--+(5pt,5pt)--+(0pt,5pt)--cycle ;
\draw (5pt,0)--+(5pt,0)--+(5pt,5pt)--+(0pt,5pt)--cycle;
\end{tikzpicture}}
$, $\Xi = \Upsilon_{d_{\iota}}$ and $\Xi' = \Upsilon_{d_{\iota'}}$.

According to what we recollected, the space~$J^{\pd}$ is spanned by $\sum_{i=1}^n \pdi{e_i}{e_i}$ on the one hand. On the other hand, the violation-free vectors are the $(\pdi{e_i}{e_j})_{1 \le i<j \le n}$, which are already orthogonal to $J^\pd$, and the $(\pdi{e_i}{ e_i})_{2 \le i \le n}$. The projection on~$\Sc^{[\pd]}(\R^n)$ of these latters are the vectors $(\pdi{e_i}{e_i} - \pdi{e_1}{ e_1})_{2 \le i \le n}$, the (dimension $(n-1)\times (n-1)$) Gram matrix of which has for inverse 
$$\begin{pmatrix} 2 & 1 & \cdots & 1 \\ 1 & 2 & \ddots & \vdots  \\ \vdots & \ddots & \ddots & 1 \\ 1 & \cdots & 1 & 2 \end{pmatrix} ^{-1} = \begin{pmatrix} \frac{n-1}{n} & -\frac{1}{n} & \cdots & -\frac{1}{n} \\ -\frac{1}{n} & \frac{n-1}{n} & \ddots & \vdots  \\ \vdots & \ddots & \ddots & -\frac{1}{n} \\ -\frac{1}{n} & \cdots & -\frac{1}{n}& \frac{n-1}{n} \end{pmatrix}. $$
In full generality $Z_{\pdi{\eps_j}{\eps_j},\pdi{\eps_{j'}}{\eps_{j'}}} (X_\De, X_{\De'})$ is equal to :
$$\sum_{2 \le i \le n}  (x_{i,j}^2-x_{1,j}^2)(x_{i,j'}'^2-x_{1,j'}'^2) - \frac{1}{n } \sum_{2 \le i , i' \le n} (x_{i,j}^2-x_{1,j}^2) (x_{i',j'}'^2-x_{1,j'}'^2)\qquad \qquad$$
$$  \qquad \qquad \qquad \qquad \qquad \qquad \qquad \qquad \qquad \qquad + \sum_{1 \le i \neq i' \le n} x_{i,j}x_{i',j} \, x_{i,j'}' x_{i',j'}'  $$
Since the orthogonal group~$\Or(n)$ acts transitively on the set of all flag, we can content ourselves in computing  $Z_{\Upsilon_{d_{\iota}},\Upsilon_{d_{\iota'}}} (\De_0, \De)$ where $\De_0$ is the flag of shape~$\la$ spanned by the vectors $(e_i)_{1 \le i \le s}$. This amounts to substitute the $x_{i,j}$ by the Kronecker symbol $\de_{i,j}$. At the end we get
$$Z_{\Upsilon_{d_{\iota}},\Upsilon_{d_{\iota'}}} (\De_0, \De) = \sum_{\substack{1 \le i \le d_\iota \\ 1 \le j \le d_{\iota'}}} x_{i,j}^2 -\frac{d_\iota}{n} \sum_{\substack{d_\iota < i \le n \\ 1 \le j \le d_{\iota'}}} x_{i,j}^2 $$
The first sum identifies with the sum of the square norms of the projection of the the $d_{\iota'}$ first vectors of $X_\De$ on $(\De_0)_\iota$ the $\iota$\textsuperscript{th} space of $\De_0$, in other words with the trace $\tr(\text{pr}_{(\De_0)_\iota} \circ \text{pr}_{(\De)_{\iota'}})$ or equivalently with the sum of the square cosines of the principle angles defined between $(\De_0)_{\iota}$ and $\De_{\iota'}$. The second sum equals~$d_{\iota'}$ since~$X_\De$ is a matrix with orthonormal columns.

In general we have thus 
\be  Z_{ \Upsilon_{d_{\iota}}, \Upsilon_{d_{\iota'}} } = \tr( \text{pr}_{\De_\iota} \circ \text{pr}_{\De_{\iota'}} ) - \frac{d_\iota d_{\iota'}}{n} . \label{eq:fz2}\ee

\section{Vexillar designs}
\label{s:desflag}

\begin{defi} We recall that the measure is normalised to one. We call \emph{vexillar $t$-design} any finite subset $\mathcal{D}$ of the flag manifold~$\vdra$ such that for any polynomial function $f$ of degree less than  or equal to~$t$, that is a function belonging to $\ds \bigoplus_{\substack{|\mu| \le t \\  \lfloor \mu \rfloor \le m }} \Sc^{[\mu]} (\R^n) ^{\oplus n_\mu} $, the quadrature formula
\be \label{e:des} \int_{\vdra} f(\De) \ud \De = \frac{1}{|\mathcal{D}|} \sum_{\De \in \mathcal{D}} f(\De) \ee
is actually a real equality. 
\end{defi}

\begin{thm} \label{thm:eqtdes} The following conditions are equivalent: 
\begin{enumerate}
\item The set $\mathcal{D}$ is a $t$-design.
\item For any homogeneous polynomial of degree less or equal than $t$, for any orthogonal transformation $\tau \in \Or(n)$, we have 
$$\sum_{\De \in \mathcal{D} } f(\De) = \sum_{\De \in \mathcal{D} }  f(\tau . \De). $$
\item For any partition~$\mu$ of degree $0<|\mu| \le t$, for any $\Xi \in  N^\mu =  \Sc^{\mu}(\R^m)^{\Or(\m)} $ and for any function $f$ of $\phi( \Sc^{[\mu]} \otimes \Xi )$, the sum $\sum_{\De \in \mathcal{D}} f(\De) $ is zero.
\item For any partition~$\mu$ of degree $0<|\mu| \le t$, for any $\Xi$ and $\Xi' \in N^\mu$ and for any $\De' \in \vdra$,
$$\sum_{\De \in \mathcal{D}} Z_{\Xi, \Xi'}(\De, \De') = 0$$
\end{enumerate}
\end{thm}

\begin{dem} We organise the proof as follows.
$$
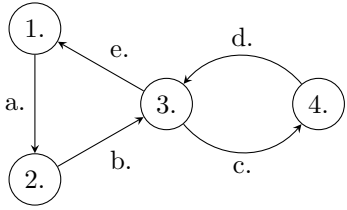
\begin{tikzpicture}
\node[draw, circle] (a) at (150:2) {$1.$};
\node[draw, circle] (b) at (-150:2) {$2.$};
\node[draw, circle] (c) at (0,0) {$3.$};
\node[draw, circle] (d) at (2,0) {$4.$}; 
\draw[-stealth] (a) -- node [anchor = east] {a.} (b);
\draw[-stealth] (b) --node [anchor = north west] {b.} (c);
\draw[-stealth] (c) .. controls +(-50:1) and +(-130:1) .. node [anchor = north] {c.} (d);
\draw[-stealth] (d) .. controls +(130:1) and +(50:1) .. node [anchor = south] {d.} (c);
\draw[-stealth] (c) -- node [anchor = south west] {e.}(a);
\end{tikzpicture} $$

\begin{itemize}
\item[a.]  Since the measure $\ud \De$ is invariant under the variable change $\De \mapsto \tau. \De$, we have, using equation (\ref{e:des}) twice,
$$\sum_{\De \in \mathcal{D} } f(\De) = |\mathcal{D} |\int_{\vdra} f(\De) \ud \De = |\mathcal{D} | \int_{\vdra} f(\tau . \De) \ud \De = \sum_{\De \in \mathcal{D} }  f(\tau . \De) $$
\item[b.] The map $\phi( \Sc^{[\mu]} \otimes \Xi ) \to \R$ defined by $f \mapsto \sum_{\De \in \mathcal{D}} f(\De)$ is a linear application. Its kernel is an  $\Or(n)$-irreducible subspace which is non trivial because of the rank theorem. It can be but the whole space $\phi( \Sc^{[\mu]} \otimes \Xi )$, which accounts for the vanishing of the sum for any~$f$.
\item[c.] This comes from $Z_{\Xi, \Xi'}(\cdot, \De')$ belonging to the space $\phi( \Sc^{[\mu]} \otimes \Xi )$.
\item[d.] The set of the functions $Z_{\Xi, \Xi'}(\cdot, \De')$ spans under the action of~$\Or(n)$ an irreducible non-zero space, which as a consequence can only be $\phi( \Sc^{[\mu]} \otimes \Xi )$. Thus the vanishing of the sum extends to the whole space.
\item[e.] Let $f$ be a function of degree less than or equal to~$t$, which we decompose by projection on the isotypic component into  $\ds f = \sum_{\deg \mu \le t} f_\mu$. The integral can be computed as follows
$$ \int_{\vdra} f(\De) \ud \De=  \psll{f}{\mathbf{1}}=  \psll{f_0}{\mathbf{1}}  = \frac{1}{|\mathcal{D}|} \sum_{\De \in \mathcal{D}} f_0(\De) = \frac{1}{|\mathcal{D}|} \sum_{\De \in \mathcal{D}} f(\De) $$
taking into account the hypothesis and the pairwise orthogonality of the spaces  $\Sc^{[\mu]} (\R^n) ^{\oplus n_\mu}$
\end{itemize}
\end{dem}
\bigskip

As in  \cite{LST} and \cite{BCN}, we have the following theorem, which connects $t$-designs with group theory. 
\begin{thm} \label{thm:eqdesgrp}Let $G$ be a finite subgroup of the orthogonal group $\Or(n)$, then the following properties are equivalent :
\begin{enumerate}
\item The decomposition of the vector space $\ds \bigoplus_{\substack{| \mu| \le t\\ \lfloor \mu \rfloor \le m} } \Sc^{[\mu]} (\R^n) ^{\oplus n_\mu} $ as invariant $G$-modules discloses the trivial representation $\bf{1_G}$ only once.
\item For any flag~$\De_0$, the orbite~$G\cdot\De_0$ of $\De_0$ under the action of $G$ forms a $t$-design.   
\end{enumerate}
\end{thm}

\begin{dem} For any partition~$\mu$ of degree $0<|\mu|\le t$ and elements~$\Xi$ and~$\Xi'$ of~$N^\mu$, we consider the application $\ds \De \mapsto \sum_{g \in G} Z_{\Xi, \Xi'}(g.\De, \De')$, where $\De'\in \vdra$ is fixed. This is a $G$-invariant map of $\phi( \Sc^{[\mu]} (\R^n) \otimes \Xi )$. But because of 1., only the constant functions can be $G$-invariant and non zero at the same time. Thus $\ds \sum_{\De \in G\cdot \De_0} Z_{\Xi, \Xi'}(\De, \De') = 0$, which is a necessary and sufficient condition (cf. 4. of theorem \ref{thm:eqtdes}) for $G\cdot \De_0$ to be $t$-design. 

Conversely, if the trivial representation $\bf{1_G}$ appears more than once among $\ds \bigoplus_{\substack{| \mu| \le t\\ \lfloor \mu \rfloor \le m}  } \Sc^{[\mu]} (\R^n) ^{\oplus n_\mu} $, it is because there is a non-zero  $G$-invariant function that belongs to some subspace $\phi( \Sc^{[\mu]} \otimes \Xi )$ with $0< | \mu |\le t$. For some flag $\De_0$, $f(\De_0)$ is non zero and $\sum_{\De \in G \cdot \De_0} f(\De) = |G|\cdot f(\De_0)$ does not vanish neither. Thus $G \cdot \De_0$ is not a $t$-design.
\end{dem}

\begin{cor} \label{cor:gradra} For the orbit~$G\cdot \De_0$ of a flag~$\De_0$ under the action of a finite subgroup~$G$ of the orthogonal group~$\Or(n)$ to constitute vexillar $t$-design, it is necessary and sufficient that the orbit of the subspace of maximal dimension~$(\De_0)_1$ be a Grassmannian $t$-design.
\end{cor}

\begin{proof} Up to multiplicities, the non empty isotypic components which appear in the decomposition of $\ds \bigoplus_{\substack{| \mu| \le t\\ \lfloor \mu \rfloor \le m} } \Sc^{[\mu]} (\R^n) ^{\oplus n_\mu} $ remain the same when we restrict to the Grassmanians of dimension $m$. \end{proof}

\begin{ex} The investigation of Grassmanian designs that arise as orbits of a sole subspace under the action of a group has been undertaken in different article, in particular~\cite{BCN} and especially in~\cite{B} for a general overview in lower dimension. 
\end{ex}

\section{Strongly perfect lattices}
\label{s:sperf}

We fix in this section a partition $\la$ with less than $n$ parts. We suppose that the parts $(\breve{\la}_i)_{1 \le i \le \breve{s}}$ of $\breve{\la}$ exactly take all the values $(d_i)_{1 \le i \le \ell}$ but can possibly be repeated. The determinant of a lattice is the determinant of the Gram matrix of one of its basis. To any lattice $L$ contained in $\R^n$, we recall that the \emph{Hermite invariant} $\ga(L)$ is given by 
\be \ga(L) = \inf_{\La \subset L} \frac{\det(\La_1) \cdots \det(\La_{\breve{s}}) }{\left( \det (L) \right) ^{\frac{|\la|}{n}}}\ee
where $\La = (\La_1, \dots, \La_{\breve{\la}})$ is a chain of nested sub-lattices of $L$ satisfying the condition $\text{rk} (\La_i) = \breve{\la}_i$ for all $1 \le i \le \breve{s}$. In the setting of \cite{W}, the partition~$\la$ is actually the weight of a representation of the general linear group.

We call \emph{minimal flag} and we denote~$\sL$ their set any chain of nested lattices  that achieves the minimum of $\ga(L)$. The set $\sL$ is finite ; its cardinal is denoted $\sD$.

With a flag of lattices $\La = (\La_1, \dots, \La_{\breve{\la}})$, we associate the \emph{sum of the orthogonal projections} $\textrm{pr}_{\La_i}$ on the vector space spanned by $\La_i$,
$$ \Pi_\La = \sum_{i=1}^{\breve{s}} \textrm{pr}_{\La_i},$$
which is a symmetric endomorphism.

\begin{defi} \begin{enumerate} \item A lattice is called \emph{perfect} with respect to $\la$ if the endomorphisms $(\Pi_\La )_{\La \in \sL } $ span the space of symmetric endomorphisms.
\item A lattice is called \emph{eutactic} with respect to $\la$ if the identity is a linear combination with only positive coefficients of all the sums of projections $(\Pi_\La )_{\La \in \sL } $.
\item A lattice is called \emph{extreme} with respect to $\la$ if it achieves a local maximum of $\ga(L)$.
\end{enumerate}
\end{defi}

The relevance of these definition is rooted in the following \`a la Voronoi result.
\begin{thm} (\cite{Mey}) \label{thm:vor} In order that a lattice be extreme with respect to $\la$, it is necessary and sufficient that it be perfect and eutactic with respect to $\la$.
\end{thm}

\begin{defi} \begin{enumerate} \item We say that a lattice is \emph{strongly perfect} with respect to $\la$ if the set of all its minimal flags of shape $\la$ carries a $4$-design.
\item We say that a lattice is \emph{strongly eutactic} with respect to $\la$ if it is eutactic and all the eutaxy coefficients are identical. 
\end{enumerate}
\end{defi}

\begin{prop} \label{prop:feut} A lattice the minimal flags of which form a $2$-design is strongly eutactic with respect to $\la$.
\end{prop}

\begin{dem} We shall show that the sum $ S = \sum_{\La \in \sL } \Pi_\La$ is proportionnal to the identity. Precisely, a trace computation points to the sole possible constant $c$ being equal to $\frac{|\la|\cdot  \sD}{n}$. Since the collection of the sums of projections~$\Pi_{\De_0}$ spans the space of symmetric endomorphisms when $\De_0$ runs through the set of flags $\vdra$, it suffices to check that $\psll{ S}{\Pi_{\De_0} } = \frac{\sD \cdot |\la|^2}{n}$. Since on the one hand 
$$  \psll{\Pi_\La}{\Pi_{\De_0}} = \sum_{1 \le i,j \le s} \left( Z_{\Upsilon_{\breve{\la}_i},\Upsilon_{\breve{\la}_i}} (\La,\De_0)+ \frac{ \breve{\la}_j \breve{\la}_j}{n}  \right)$$
according to the calculation \ref{eq:fz2} and on the other hand according to the theorem~\ref{thm:eqtdes} of caracterisation of designs the sum $ \sum_{\La \in \sL } Z_{\Upsilon_{\breve{\la}_i},\Upsilon_{\breve{\la}_i}} (\La,\De_0)$ vanishes, this equality is fulfilled.
\end{dem}

\begin{thm} A strongly perfect lattice with respect to $\la$ is extreme with respect to the partition $\la$.
\end{thm}

\begin{dem} Taking into account theorem \ref{thm:vor} and proposition \ref{prop:feut}, we are still left to prove that strong perfection implies perfection, or in other words that minimal flags forming a $4$-design implies the rank of the applications $(\Pi_\La)_{\La \in \sL}$ being maximal.

To that end, we show that the matrix $C = \big( \psll{ \Pi_\De}{ \Pi_{\De'}}  \big)_{\De,\De' \in \sL }$ has rank $n \, (n+1)  /2$. Let us denote in a first time $P$ the matrix $P_{\Upsilon_{\breve{\la}_{\iota_1}}, \Upsilon_{\breve{\la}_{\iota_2}}}  = \big(   Z_{\Upsilon_{\breve{\la}_{\iota_1}}, \Upsilon_{\breve{\la}_{\iota_2}}} (\De, \De')  \big)_{\De,\De' \in \sL }$ and compute the matrix product $P_{\Upsilon_{\breve{\la}_{\iota_1}}, \Upsilon_{\breve{\la}_{\iota_2}}}  \cdot P_{\Upsilon_{\breve{\la}_{\iota_3}}, \Upsilon_{\breve{\la}_{\iota_4}}} $
\Bea
&& \sum_{\De'' \in \sL}  Z_{\Upsilon_{\breve{\la}_{\iota_1}}, \Upsilon_{\breve{\la}_{\iota_2}}} (\De, \De'') Z_{\Upsilon_{\breve{\la}_{\iota_3}}, \Upsilon_{\breve{\la}_{\iota_4}}} (\De'', \De') \\
&=&
\frac{1}{N^2} \sum_{\De'' \in \sL}  \sum_{1\le i,j \le N} (h_i \otimes \Upsilon_{\breve{\la}_{\iota_1}})(\De)  (h_i \otimes \Upsilon_{\breve{\la}_{\iota_2}})(\De'')  (h_j \otimes \Upsilon_{\breve{\la}_{\iota_3}})(\De'')  (h_j \otimes \Upsilon_{\breve{\la}_{\iota_4}})(\De')\\
&=& \frac{1}{N^2} \sum_{1 \le i,j \le N}  (h_i \otimes \Upsilon_{\breve{\la}_{\iota_1}})(\De)  (h_j \otimes \Upsilon_{\breve{\la}_{\iota_4}})(\De') \sum_{\De'' \in \sL }  (h_i \otimes \Upsilon_{\breve{\la}_{\iota_2}})(\De'')  (h_j \otimes \Upsilon_{\breve{\la}_{\iota_3}})(\De'').\\ 
\Eea
However, the application $f : \De'' \mapsto (h_i \otimes \Upsilon_{\breve{\la}_{\iota_2}})(\De'')  (h_j \otimes \Upsilon_{\breve{\la}_{\iota_3}})(\De'')$ is of degree 4 and can decompose into a sum on its isotypic components
$$f = f_0 + f_\pd + f_\pdd + f_\pq .$$
Since $\mathcal{S}(\De)$ is a $4$-design, the sum $\sum_{\De'' \in \sL }  f_\pd(\De'') + f_\pdd(\De'') + f_\pq(\De'')$ vanishes. Besides, $f_0 = \psll{ f}{ \mathbf{1} } \mathbf{1}  = \psll{ h_i \otimes \Upsilon_{\breve{\la}_{\iota_2}}}{ h_j \otimes \Upsilon_{\breve{\la}_{\iota_3}} } \mathbf{1} =  \frac{ \kappa_{i,j}}{n}  \de_{i,j} \mathbf{1} $, where $\de_{i,j}$ is the Kronecker symbol and $\kappa_{i,j} =n\cdot  \min(\breve{\la}_{i} , \breve{\la}_{j})) - \breve{\la}_{i} \breve{\la}_{j}$.
\Bea \sum_{\De'' \in \sL }  (h_i \otimes \Upsilon_{\breve{\la}_{\iota_2}})(\De'')  (h_j \otimes \Upsilon_{\breve{\la}_{\iota_3}})(\De'') 
&=& \sum_{\De'' \in \sL }  f_0(\De'')  \\
&=&  \frac{\sD \cdot  \kappa_{\iota_2, \iota_3} }{n} \de_{i,j}. \Eea
Eventually, 
$$ P_{\Upsilon_{\breve{\la}_{\iota_1}}, \Upsilon_{\breve{\la}_{\iota_2}}}  \cdot P_{\Upsilon_{\breve{\la}_{\iota_3}}, \Upsilon_{\breve{\la}_{\iota_4}}} = \frac{\sD \cdot \kappa_{\iota_2, \iota_3} }{Nn} P_{\Upsilon_{\breve{\la}_{\iota_1}}, \Upsilon_{\breve{\la}_{\iota_4}}} .$$

Then, denoting  $J$ the matrix with entries all equal to 1, 
\be  \label{eq:c} C = \sum_{1 \le i, j \le s} P_{\Upsilon_{\breve{\la}_i}, \Upsilon_{\breve{\la}_j }} + \frac{|\la|^2}{n} J \ee
and, using the relations $J  P_{\Upsilon_{\breve{\la}_{\iota_1}}, \Upsilon_{\breve{\la}_{\iota_2}}} =  P_{\Upsilon_{\breve{\la}_{\iota_1}}, \Upsilon_{\breve{\la}_{\iota_2}}} J = 0$, which stem from $\sL$ being a $2$-design, and posing  $\kappa = \sum_{1 \le i, j \le s}  \kappa_{\la_i, \la_j}$, so that 
\be \label{eq:p} \left( \sum_{1 \le i, j \le s} P_{\Upsilon_{\breve{\la}_i}, \Upsilon_{\breve{\la}_j }} \right)^2 = \frac{\sD \cdot \kappa}{Nn } \sum_{1 \le i, j \le s} P_{\Upsilon_{\breve{\la}_i}, \Upsilon_{\breve{\la}_j }}, \ee
we end up by squaring with the relation
$$ C^2 =  \frac{\sD \kappa}{Nn}   \sum_{1 \le i, j \le s} P_{\Upsilon_{\breve{\la}_i}, \Upsilon_{\breve{\la}_j }}+ \frac{|\la|^4}{n^2} \sD J =\frac{ \sD \kappa  }{Nn}C + \frac{\sD}{n^2} \left(|\la|^4 - \frac{\kappa |\la|^2}{N}\right) J   .$$
Finally, 
\be \label{eq:cc}  C^2 = \frac{\sD \kappa}{n N} C + \frac{\sD |\la|^2}{n^2} \left( |\la|^2 - \frac{\kappa}{N} \right) J. \ee
Since the matrices $C^2$, $C$ and $J$ commute and are symmetric, they can be simultaneously diagonalised. The matrix~$J$ has two eigenvalues which are~$\sD$ with multiplicity~1 and~0 with multiplicity~$\sD-1$. Let us call $(\om_i)_{1 \le i \le s}$ the eigenvalues of~$C$, where $\om_1$ is the eigenvalue associated with $(1,\dots, 1)$. Since $n CJ = |\la|^2 \sD J$ holds by multiplying (\ref{eq:c}) by $J$ we have $\om_1 = \frac{|\la|^2 \sD}{n} $. Next, using (\ref{eq:cc}), for any $i\ge 2$, $n^2 \om_i^2 = \frac{\kappa \sD  n}{N} \om_i $, whence the admissible eigenvalues are $\om_i=0$ or else $\om_i = \frac{\sD \kappa}{n  N }$. But, denoting $\al$ then multiplicity of this eigenvalue,  we obtain by a trace computation 
$$ \tr(C)= \sD \sum_{i=1}^s  (2i-1) \breve{\la}_i = \frac{|\la|^2 \sD}{n} + \al \frac{\sD \kappa}{n  N } $$
We get for $\al$ 
$$ \al = \frac{n N}{\kappa} \left( \sum_{i=1}^s  (2i-1) \breve{\la}_i -  \frac{|\la|^2}{n} \right) = N = \frac{n (n+1)}{2}- 1$$
Since the $(\kappa_{i,j})$ and $\kappa$ itself are positive, the rank of $C$ is exactly equal to $\frac{n (n+1)}{2}$ and the lattice~$L$ is perfect. \end{dem}

\section{Examples}

In each of the following cases, we apply the criterion of theorem \ref{thm:eqdesgrp} or its corollary to the automorphism group the lattices to show that any orbit is at least a 4-design. As a consequence, the lattices are strongly perfect and thus extreme.

This method applies especially well to higher dimensional lattices, where even computing the usual minimum is sometimes out of reach, or infinite famillies of lattices. We name just a few of them : some of the following lattices were already mentionned in \cite{LST} for the classical case and are actually extreme for any $\la$.

\begin{thm} \label{p:rtlat} We suppose in each case that the partition $\la$ satisfies $\lfloor \la \rfloor < n/2$.
 \begin{enumerate}
\item The root lattices $D_4$, $E_6$, $E_7$, $E_8$,are extreme with respect to any $\la$.
\item   The Thompson-Smith 248-dimensional lattice, the Fisher 78-dimensional lattice constructed by Schroeder, the 52-dimensional lattice related to the group $2 \cdot F_4(2)$ or the Leech lattice are extreme with respect to any $\la$.
\item The Barnes--Wall lattices are are extreme with respect to any $\la$
\end{enumerate}
\end{thm}

\begin{proof} \begin{enumerate} \item As it was done in \cite{BCN}, to check that the trivial representation of a group $G$ in $\ds \bigoplus_{\substack{| \mu| \le 2t\\ \lfloor \mu \rfloor \le m} } \Sc^{[\mu]} (\R^n) ^{\oplus n_\mu} $ occurs only once is equivalent to compare its number of occurence in ${\sf Sym}^t({\sf Sym}^2(\R^n))$ with the one of $\Or(n)$ since the irreducible component that are involved are the same. Character theory provides a good tool to achieve the comparison.

Under the automorphism group of $D_4$, $E_6$ and $E_8$, a flag of the corresponding space is always a 2-design. Under the automorphism group of $E_8$, a flag of $\R^8$ is a 6-design.  Under the automorphism group of the Leech lattice, a flag of $\R^{24}$ is a 10-design. At any rate, in all the cases the design is acute enough for the lattice to be strongly perfect and thus extreme. This was already checked in \cite{B} for the Grassmannian case.

\item These lattices  are related to particular groups. It was proven in \cite{TIE} that their automorphism groups afford Grassmannian 4-designs. Thus they also give rise to vexillar 4-designs and the associated lattices are extreme.

\item Let us denote $\mathcal{G}_k$ the automorphism group of the Barnes--Wall lattice of dimension $n=2^k$. According to theorem 5.1 of \cite{B}, we know that for any $k\ge 3$ and $d \le 6$, the following invariants are equal
$$\left( (\R^n)^{\otimes d} \right)^{\mathcal{G}_k} = \left( (\R^n)^{\otimes d} \right)^{\Or(n)} .$$
Since all the irreducible representations we need to take into account are included in these tensor product, we can derive that the minimal flag of the Barnes-Wall lattices are $6$-designs. As a consequence, they are strongly perfect and extreme.
\end{enumerate}
\end{proof}

\begin{prop} The lattices $K'_{10}$, ${K'_{10}}^*$ and the Coxeter-Todd lattice $K_{12}$ are strongly perfect with respect to a partition $\la$ that  has just one part but not in general.
\end{prop}

\begin{proof} The same manipulations as for theorem \ref{p:rtlat} can be performed. Yet, it turns out that there exist an invariant polynomial relative to the partition~$\pdd$. Thus the minimal vectors and there symmetric powers are $4$-designs but it is not the case for any other shape of~$\la$  (see \cite{B}).
\end{proof}

\bibliographystyle{alpha}
\bibliography{DesignDrap.bib}

\begin{thebibliography}{dCEP80}

\bibitem[Bac05]{B}
Christine Bachoc.
\newblock Designs, groups and lattices.
\newblock {\em J. Th\'eor. Nombres Bordeaux}, 17(1):25--44, 2005.

\bibitem[BBC04]{BBC}
Christine Bachoc, Eiichi Bannai, and Renaud Coulangeon.
\newblock Codes and designs in {G}rassmannian spaces.
\newblock {\em Discrete Math.}, 277(1-3):15--28, 2004.

\bibitem[BCN02]{BCN}
Christine Bachoc, Renaud Coulangeon, and Gabriele Nebe.
\newblock Designs in {G}rassmannian spaces and lattices.
\newblock {\em J. Algebraic Combin.}, 16(1):5--19, 2002.

\bibitem[BV88]{BV}
Winfried Bruns and Udo Vetter.
\newblock {\em Determinantal rings}, volume 1327 of {\em Lecture Notes in
  Mathematics}.
\newblock Springer-Verlag, Berlin, 1988.

\bibitem[dCEP80]{DCEP}
Corrado de~Concini, David Eisenbud, and Claudio Procesi.
\newblock Young diagrams and determinantal varieties.
\newblock {\em Invent. Math.}, 56(2):129--165, 1980.

\bibitem[dCP76]{DCP}
Corrado de~Concini and Claudio Procesi.
\newblock A characteristic free approach to invariant theory.
\newblock {\em Advances in Math.}, 21(3):330--354, 1976.

\bibitem[DGS77]{DGS}
Philippe Delsarte, Jean-Marie Goethals, and Johan~Jacob Seidel.
\newblock Spherical codes and designs.
\newblock {\em Geometriae Dedicata}, 6(3):363--388, 1977.

\bibitem[DRS74]{DRS}
Peter Doubilet, Gian-Carlo Rota, and Joel Stein.
\newblock On the foundations of combinatorial theory. {IX}. {C}ombinatorial
  methods in invariant theory.
\newblock {\em Studies in Appl. Math.}, 53:185--216, 1974.

\bibitem[Ful97]{Fu}
William Fulton.
\newblock {\em Young tableaux}, volume~35 of {\em London Mathematical Society
  Student Texts}.
\newblock Cambridge University Press, Cambridge, 1997.
\newblock With applications to representation theory and geometry.

\bibitem[JC74]{JC}
Alan~T. James and A.~G. Constantine.
\newblock Generalized {J}acobi polynomials as spherical functions of the
  {G}rassmann manifold.
\newblock {\em Proc. London Math. Soc. (3)}, 29:174--192, 1974.

\bibitem[LST01]{LST}
Wolfgang Lempken, Bernd Schr{\"o}der, and Pham~Huu Tiep.
\newblock Symmetric squares, spherical designs, and lattice minima.
\newblock {\em J. Algebra}, 240(1):185--208, 2001.
\newblock With an appendix by Christine Bachoc and Tiep.

\bibitem[LT85]{LT}
Glenn Lancaster and Jacob Towber.
\newblock Representation-functors and flag-algebras for the classical groups.
  {II}.
\newblock {\em J. Algebra}, 94(2):265--316, 1985.

\bibitem[Mey09]{Mey}
Bertrand Meyer.
\newblock Generalised {H}ermite constants, {V}oronoï theory and heights on flag
  varieties.
\newblock {\em Bull. Soc. Math. Fr.}, 137(1), 2009.

\bibitem[Tie06]{TIE}
Pham~Huu Tiep.
\newblock Finite groups admitting {G}rassmannian 4-designs.
\newblock {\em J. Algebra}, 306(1):227--243, 2006.

\bibitem[Ven01]{V}
Boris Venkov.
\newblock R\'eseaux et designs sph\'eriques.
\newblock In {\em R\'eseaux euclidiens, designs sph\'eriques et formes
  modulaires}, volume~37 of {\em Monogr. Enseign. Math.}, pages 10--86.
  Enseignement Math., Geneva, 2001.

\bibitem[Wat00]{W}
Takao Watanabe.
\newblock On an analog of {H}ermite's constant.
\newblock {\em J. Lie Theory}, 10(1):33--52, 2000.

\end{thebibliography}

\end{document}